\documentclass[11pt,DIV11]{article}
\usepackage[utf8]{inputenc}
\usepackage[T1]{fontenc}
\usepackage{lmodern}

\usepackage{amssymb,amsmath,amsthm}
	\newtheorem{example}{Example}[section]
	\newtheorem{theorem}{Theorem}[section]
	\newtheorem{lemma}{Lemma}[section]
	\newtheorem{proposition}{Proposition}[section]
	\newtheorem{definition}{Definition}[section]
	
		\newtheorem{assumption}{Assumption}[section]
	\numberwithin{equation}{section}
	
\usepackage{xcolor}
\usepackage{dsfont}

\usepackage{graphicx}

\usepackage{mathrsfs} 

\renewcommand{\sc}{\textsc} 

\newcommand{\R}{{\mathbb R}}				
\newcommand{\N}{{\mathbb N}}				
\newcommand{\st}{\text{s.t. }}			
\newcommand{\abs}[1]{\left|#1\right|}	

\newcommand{\tspace}{[0,1]}					
\newcommand{\yspace}{{\R^n}}				
\newcommand{\uspace}{{\R^m}}				

\newcommand{\ringu}{\mathcal{U}}			

\newcommand{\cuspace}{{\beta_\ringu\uspace}}	

\newcommand{\tm}{\tau}				
\newcommand{\um}{\omega}				

\usepackage[normalem]{ulem}
\newcommand{\new}[1]{{\color{black}#1}}

\author{Didier Henrion\footnote{LAAS-CNRS, Universit\'e de Toulouse, CNRS, Toulouse, France \& Faculty of Electrical Engineering, Czech Technical University in Prague, Prague, Czechia (\texttt{henrion@laas.fr})},\;
Martin Kru\v{z}\'{\i}k\footnote{Czech Academy of Sciences, Institute of Information Theory and Automation, Praha, Czechia \& Faculty of Civil Engineering, Czech Technical University in Prague, Prague, Czechia (\texttt{kruzik@utia.cas.cz})},\;
Tillmann Weisser\footnote{LAAS-CNRS, Universit\'e de Toulouse, CNRS, Toulouse, France (\texttt{tweisser@laas.fr})}}

\title{Optimal control problems with oscillations, concentrations and discontinuities}

\begin{document}
\maketitle
\begin{abstract}
Optimal control problems with oscillations (chattering controls) and concentrations (impulsive controls) can have integral performance criteria such that concentration of the control signal occurs at a discontinuity of the state signal. Techniques from functional analysis (anisotropic parametrized measures) are applied to give a precise meaning of the integral cost and to allow for the sound application of numerical methods. We show how this can be combined with the Lasserre hierarchy of semidefinite programming relaxations.
\end{abstract}

\section{Introduction}
As a consequence of optimality, various limit behaviours can be observed in optimal control: minimizing control law sequences may feature increasingly fast variations, called oscillations (chattering controls \cite{Young1969}), or increasingly large values, called concentrations (impulsive controls \cite{Luenberger1969}). The simultaneous presence of oscillations and concentrations in optimal control needs careful analysis and specific mathematical tools, so that the numerical methods behave correctly. Previous work of two of the authors \cite{Claeys2017} combined tools from partial differential equation analysis (DiPerna-Majda measures \cite{DiPerna1987}) and semidefinite programming relaxations (the moment-sums-of-squares or Lasserre hierarchy \cite{Lasserre2008}) to describe a sound numerical approach to optimal control in the simultaneous presence of oscillations and concentrations. To overcome difficulties in the analysis, a certain number of technical assumptions were made, see \cite[Assumption 1, Section 2.2]{Claeys2017}, so as to avoid the simultaneous presence of concentrations (in the control signals) and discontinuities (in the system trajectories).

In the present contribution we remove these technical assumptions and accommodate the simultaneous presence of concentrations and discontinuities, while allowing oscillations as well. For this, we exploit a recent extension of the notion of DiPerna-Majda measures called anisotropic parametrized measures \cite{Kalamajska2017}, so that it makes sense mathematically while allowing for an efficient numerical implementation with semidefinite programming relaxations.

To motivate further our work, let us use an elementary example to illustrate the difficulties that may be faced in the presence of discontinuities and concentrations. Consider the optimal control problem
\begin{equation}\label{eq:ipp}
\begin{split}
\inf_u &\int_0^1  (t+y(t))u(t) dt\\
\st & \dot{y}(t)= u(t),\quad y(0)=0,\quad y(1)=1,\\
& 1 \geq y(t) \geq 0, \quad u(t) \geq 0, \quad t \in [0,1]\\
\end{split}
\end{equation}
where the infimum is with respect to measurable controls of time.
The trajectory $y$ should move the state from zero at initial time to one at final time, yet for the non-negative integrand to be as small as possible, the control $u$ should be zero all the time, except maybe at time zero. We can design a sequence of increasingly large controls $u$ that drive $y$ from zero to one increasingly fast. We observe that this sequence has no limit in the space of measurable functions but it tends (in a suitable weak sense) to the Dirac measure at time zero. We speak of control signal concentration or impulsive control. The integrand contains the product $yu$ of a function whose limit becomes discontinuous at a point where the other function has no limit, hence requiring careful analysis. Here however, this product can be written $y\dot{y}=\tfrac{d}{dt}\tfrac{y^2}{2}$ and hence the integral term is well defined since $\int_0^1 y\dot{y}dt = \frac{y(1)^2-y(0)^2}{2} = \frac{1}{2}$. Consequently the cost in \eqref{eq:ipp} is equal to $\int_0^1 t u(t) dt+ \tfrac{1}{2}$ and independent of the actual trajectory.

This reasoning is valid because $\dot{y}(t)=u(t)$ in problem (\ref{eq:ipp}), but this integration trick cannot be carried out for more general differential equations. For example we cannot solve analytically the following modified optimal control problem
\begin{equation}\label{ex}
\begin{split}
\inf_u &\int_0^1  (t+y(t))u(t) dt\\
\st & \dot{y}(t)= \sqrt{\varepsilon^2+u^2(t)},\quad y(0)=0,\quad y(1)=1,\\
& 1 \geq y(t) \geq 0, \quad u(t) \geq 0, \quad t \in [0,1]\\
\end{split}
\end{equation}
where $\varepsilon$ is a given real number.
Providing a mathematically sound framework for the analysis of this kind of phenomenon combining concentration and discontinuity, and possibly also oscillation (not illustrated by the simple example above), is precisely the purpose of our paper.

\section*{\new{Contribution}}

\new{
The contribution of our paper with respect to previous work can be summarized as follows:
\begin{itemize}
\item we propose a unified approach for handling the simulatenous presence of oscillations, concentrations and discontinuities, where previous work considered either oscillations without concentrations (see  \cite{Vinter1993,Roubicek1997,Fattorini99,Gaitsgory2009} and references therein), or concentrations without oscillations (see \cite{Claeys2014} and references therein), or oscillations and concentrations without discontinuities (see \cite{Claeys2017} and references therein);
\item we remove the technical assumptions of \cite{Claeys2017} to allow for the simultaneous presence of concentration (of the control) and discontinuity (of the trajectory);
\item as in \cite{Claeys2017}, our approach allows for a constructive solution  via the Lasserre hierarchy \cite{Lasserre2008}; this now provides a unified numerical scheme to deal with oscillations, concentrations and discontinuities;
\item we make a connection between anisotropic measures and the occupation measures, which are classical objects in dynamical systems and Markov decision processes, and which have been used in linear reformulations of nonlinear optimal control problems \cite{Vinter1993,Lasserre2008,Gaitsgory2009}; the notion of occupation measure was extended in \cite{Claeys2013,Claeys2014} to cope with concentration (also called implusive controls);  it was pointed out in \cite{Zidani2013} that this extension allows for optimization over all possible graph completions, a tool introduced in \cite{Bressan1988} -- see also \cite{Bressan2007} -- to deal with differential equations with discontinuous solutions. Anisotropic measures allow for a further generalization of these approaches.
\end{itemize}
}

\section*{\new{Outline}}

\new{
The outline of the paper is as follows.
In Section \ref{sec:relaxing} we describe the limit phenomena typical of optimal control, namely oscillations, concentrations and discontinuities, as well as the linear formulation of optimal control problems using measures.
In Section \ref{sec:extDiPeMa} we introduce the anisotropic parametrized measures, illustrating their use with elementary examples. We show how these measures can cope with concentrations and discontinuities, giving a meaning to otherwise ill-defined integrals.
In Section \ref{sec:OCD} we apply the anisotropic parametrized measures to optimal control, and
in Section \ref{sec:relatcontrol} we describe their relationship with occupation measures, a classical tool in dynamical systems and Markov decision processes.
In Section  \ref{sec:approach} we describe how the Lasserre hierarchy can be applied to our problem, and in Section \ref{sec:example} we provide a simple illustrative example that can be solved numerically, and then analytically.
Finally, concluding remarks are gathered in Section \ref{sec:conclu}.
}
\section{Relaxing Optimal Control}\label{sec:relaxing}
Let $L:\tspace\times\yspace\times\uspace \to \R$ and $F:\tspace\times\yspace\times\uspace \to\yspace$ be continuous functions. For initial $y_0$ and final conditions $y_1$ in $\yspace$ and some integer $1\leq p \leq \infty$, the formulation of the classical optimal control problem is 
\begin{equation}\label{eq:classical ocp}
\begin{split}
v^* := \inf_u\int_0^1 & L(t,y(t),u(t))dt \\
\st & \dot{y}(t) = F(t,y(t),u(t)),\;\; y(0)=y_0,\;\; y(1)=y_1, \\
&  y\in {\mathscr W}^{1,1}(\tspace;\yspace), \: u\in {\mathscr L}^p(\tspace;\uspace)
\end{split}
\end{equation}
where  ${\mathscr W}^{1,p}([0,1];X)$ is the space of functions from $[0,1]$ to $X$ whose weak derivative belongs to 
${\mathscr L}^p([0,1];X)$, the space of functions from $[0,1]$ to $X$ whose $p$-th power is Lebesgue integrable.

\new{A pair $(u,y)$ with a control $u \in {\mathscr L}^p(\tspace;\uspace)$ and the corresponding trajectory $y \in {\mathscr W}^{1,1}(\tspace;\yspace)$ satisfying the differential equation of problem \eqref{eq:classical ocp} is called admissible.  Given a minimizing admissible sequence  $(u_k, y_k)_{k\in\N}$,} the infimum in \eqref{eq:classical ocp} might not be attained because $(u_k)_{k\in\N}$ might not converge in ${\mathscr L}^p$ and $(\new{y_k})_{k\in\N}$ might not converge in ${\mathscr W}^{1,1}$ \new{as ${\mathscr L}^1$ is not reflexive.} To overcome this issue, it has been proposed to relax the regularity assumptions on $u$ \new{and $y$}. We discuss some of the approaches now in detail.

\subsection{Oscillations} \label{sec:oscillation}
The limit of a minimizing sequence for \eqref{eq:classical ocp} might fall out of the feasible space because of oscillation effects of  $(u_k)_{k\in\N}$. Consider for example the optimal control problem 
\begin{equation}\label{ocp:oscillation}
\begin{split}
\inf_u \int_0^1& (u(t)^2-1)^2 + y(t)^2 dt\\
\st & \dot{y}(t) = u(t),\; y(0)= 0,\; y(1)= 0, \\	
	& y\in {\mathscr W}^{1,4}([0,1]),\quad u\in {\mathscr L}^4([0,1]).
\end{split}
\end{equation}
As the integrand in the cost is a sum of squares, the value is at least zero. To see that actually it is equal to zero, consider the sequence of controls $(u_k)_{k\in\N}\subseteq {\mathscr L}^4([0,1])$ defined by 
\begin{equation}\label{eq:oscillating sequence}
u_k(t):=\left\lbrace\begin{array}{rl}
1, & \text{if } t\in \left[\tfrac{2l+1}{2^k},\tfrac{l+1}{2^{k-1}}\right],\; 0 \leq l \leq k-1\\
-1,&\text{otherwise}
\end{array} \right.
\end{equation}
for $k>1$ and $u_1:=0$. For the corresponding sequence of trajectories $(y_k)_{k\in\N}$ defined by $y_k(t) := \int_0^t u_k(s)ds$ it holds that $y_k\in {\mathscr W}^{1,4}([0,1])$ and $y_k(1)=0$ as desired. Hence, $(u_k)_{k\in\N}$ is a sequence of feasible controls. A short calculation shows that using this sequence the cost in \eqref{ocp:oscillation} converges to zero.
\new{The sequence $(y_k)_{k\in\N}$ converges to $y_\infty:=0$ in ${\mathscr W}^{1,4}$, but the sequence  $(u_k)_{k\in\N}$ does not converge to $u_\infty:=0$ in ${\mathscr L}^4$.}

In contrast to that, the sequence of measures defined by $d\nu_k(t,u):=\delta_{u(t)}(du|t)dt$ converges weakly to $d\nu(t,u):=\tfrac{1}{2}(\delta_{-1}+\delta_{1})(du)dt$ in the sense that for all $f\in {\mathscr C}(\tspace)$ and $g\in {\mathscr C}_p(\R)$: 
\begin{equation}\label{eq: convergence Young}
\lim_{k\to\infty}\int_0^1 \int_{\R} f(t)g(u)d\nu_k(t,u) = \int_0^1 \int_{\R} f(t)g(u)d\nu(t,u)
\end{equation}
where ${\mathscr C}_p(\R) := \{ g\in {\mathscr C}(\R) : g(u) = o(\abs{u}^p) \text{ for } \abs{u}\to\infty\}$ is the set of continuous functions of less than $p$-th growth. 
Integration then yields $y_\infty(1)=
\int_0^1 \int_\R u d\nu(t,u) = \int_0^1 \int_\R u \tfrac{1}{2}(\delta_{-1}+\delta_{1})(du)dt = 0$.
A similar reasoning shows that the cost with respect to $\nu$ is zero.
 
More generally, this observation motivates to relax the regularity assumptions on the control $u$  in \eqref{eq:classical ocp} and also allow for limits $d\nu(t,u) = d\um(u|t)dt$ of control sequences $(u_k)_{k\in\N}\subseteq {\mathscr L}^p(\tspace;\uspace)$. In general the measure $\um$ depends on time, i.e., we have a family of probability measures $\um(.|t)_{t\in\tspace} \subset {\mathscr P}(\R^m)$, where ${\mathscr P}(X)$ denotes the set of probability measures on $X$, i.e. non-negative Borel regular measures with unit mass. Such parametrized measures obtained as limits of a sequence of functions $(u_k)_{k\in\N}\subseteq {\mathscr L}^p(\tspace;\uspace)$ have been called $L^p$-{Young measures}. 
For an explicit characterization of these measures see e.g. \cite{Kruzik1996}. 
For a comprehensive reference on Young measures and their use in the control of ordinary and partial differential equations, see \cite[Part III]{Fattorini99}.

The relaxed version of \eqref{eq:classical ocp} that now takes into account oscillating control sequences can be written as
\begin{equation}
\begin{split}
\inf_\um \int_0^1 \int_\uspace & L(t,y(t),u)\,\um(du|t)dt  \\
\st & \int_0^1 \int_\uspace F(t,y(t),u)\,\um(du|t)dt = y_1 - y_0\\
	& y\in {\mathscr W}^{1,1}(\tspace;\yspace),\; \um(.|t)\in {\mathscr P}(\uspace)
\end{split}
\end{equation}
where the constraint is a reformulation of the differential equation
$\dot{y}(t) = \int_{\R^m} F(t,y(t),u)\omega(du|t), \: t\in\tspace$
with the boundary conditions $y(0)=y_0$ and $y(1)=y_1$.

\subsection{Concentrations} \label{sec:concentration}
Oscillation of the control sequence \new{due to nonconvexity of the functional} is not the only reason that prevents the infimum in \eqref{eq:classical ocp} of being attained. As a second example consider the following problem of optimal control:
\begin{equation}\label{ocp:concentration}
\begin{split}
\inf_u \int_0^1& (t-\tfrac{1}{2})^2u(t)dt\\
\st & \dot{y}(t) = u(t) \geq 0,\; y(0)= 0,\; y(1)= 1, \\	
	& y\in {\mathscr W}^{1,1}([0,1]), \quad u\in {\mathscr L}^1([0,1]).
\end{split}
\end{equation}
Note that the control enters into the problem linearly. The value is zero as the integrand is positive and using the sequence of controls
\begin{equation}\label{eq:concentrating sequence}
u_k(t):=\left\lbrace\begin{array}{rl}
k, & \text{if } t\in \left[\tfrac{k-1}{2k},\tfrac{k+1}{2k}\right]\\
0, &\text{else}
\end{array} \right.
\end{equation}
the cost converges to zero. Neither $(u_k)_{k\in\N}$ nor any subsequence converges in ${\mathscr L}^1$ \new{as this space is not reflexive.} In contrast to the previous example this time $(y_k)_{k\in\N}$ does not converge in ${\mathscr W}^{1,1}([0,1])$ neither \new{because ${\mathscr W}^{1,1}$ is not reflexive}. We hence use the extension $\mathscr{BV}([0,1])$, the space of functions with bounded variation, as a relaxed space for the trajectory. Following the same approach as before we consider the control as a measure $d\nu_k(t,u):=\delta_{u_k(t)}(du)dt$. As $u$ appears linearly in \eqref{ocp:concentration} we can directly integrate with respect to $u$ and define a sequence of probability measures $(\tm_k)_{k\in \N}\subseteq{\mathscr P}(\tspace)$ by $\tm_k(dt) := \int_{\R} ud\nu_k(t,u)$. A short calculation shows that this sequence has the weak limit $\tm:=\delta_{\tfrac{1}{2}}$, i.e. for all $f\in {\mathscr C}([0,1])$ it holds
$\lim_{k\to\infty} \int_0^1 f(t)\tm_k(dt) = \int_0^1 f(t) \tm(dt)$.
Note that by integrating before passing to the limit we transfer the unboundedness of the control into the measurement of time and only keep the direction (i.e. $+1$ in this example) of the control. 
Whereas we observed a superposition of two different controls in the previous example, here we see a concentration of the control in time. For optimal control problems with linear growth in the control:
\[
\begin{split}
\inf_u\int_0^1 & L(t,y(t))u(t)dt \\
\st & \dot{y}(t) = F(t,y(t))u(t),\;\; y(0)=y_0,\;\; y(1)=y_1, \\
& y\in {\mathscr W}^{1,1}(\tspace;\yspace), \quad u\in {\mathscr L}^1(\tspace;\uspace)
\end{split}
\]
we can therefore build the following relaxation that can take into account concentration effects of the control:
\begin{equation}
\begin{split}
\inf_\tm \int_0^1 & L(t,y(t))\tm(dt) \\
\st & \int_0^1 F(t,y(t))\tm(dt) = y_1 - y_0, \\
    & y\in \mathscr{BV}(\tspace;\yspace), \quad\tm\in{\mathscr P}(\tspace). 
\end{split}
\end{equation}
See \cite{Claeys2014} for an application of the moment-sums-of-squares hierarchy for solving numerically non-linear control problems in the presence of concentration.

\subsection{Oscillation and Concentration} \label{sec:oscillation and concentration}
The relaxations proposed so far allow to consider controls that are either oscillating in value or concentrating in time. However it is possible that both effects appear in the same problem. Consider for example
\begin{equation}\label{ocp:oscillation and concentration}
\begin{split}
\inf_u \int_0^1& \frac{u(t)^2}{1+u(t)^4}+\left(y(t)-t\right)^2 dt\\
\st & \dot{y}(t) = u(t) \geq 0,\; y(0)= 0,\; y(1)= 1, \\	
	& y\in {\mathscr W}^{1,1}([0,1]),\quad u\in {\mathscr L}^1([0,1]).
\end{split}
\end{equation}
The infimum value zero of \eqref{ocp:oscillation and concentration} can be approached arbitrarily close by a sequence of controls $(u_k)_{k\in\N}$ defined by
\begin{equation}\label{eq:oscillating and concentrating sequence}
u_k(t):=\left\lbrace\begin{array}{rl}
k, & \text{if } t\in \left[\tfrac{l}{k}-\tfrac{1}{2k^2},\tfrac{l}{k}+\tfrac{1}{2k^2}\right],\;1\leq l< k\\
0, &\text{else}
\end{array} \right.
\end{equation}
for $k>1$ and $u_1:=1$. The idea to capture the limit behaviour of this sequence is to combine a Young measure on the control and replacing the uniform measure on time by a more general measure on time. Note that due to linearity it was possible in Section \ref{sec:concentration} to transfer the limit behaviour of the control into the measurement of time. In the present example the control enters non-linearly in the cost, which is why we will need to allow the control to take values at infinity. We consider a metrizable compactification $\beta_\ringu\R$ of the control space corresponding to the ring $\ringu$ of complete and separable continuous functions (see Section \ref{sec:DiPerna Majda measures} for more details). Then the sequence of measures $d\nu_k(t,u):=\delta_{u_k(t)}(du|t)dt$ converges to $d\nu(t,u) := \um(du)\tm(dt)$ with $\um(du):=\tfrac{1}{2}(\delta_{0}+\delta_\infty)(du)$ and $\tm(dt):=2dt$ understood in the following weak sense for all 
$f\in {\mathscr C}(\tspace)$ and $g_0\in\ringu$: 
\begin{equation}\label{eq: convergence DPM}
\begin{split}
\lim_{k\to\infty}\int_0^1 \int_\R f(t)g_0(u)(1+\abs{u}^p) d\nu_k(t,u) =  \\
\int_0^1 \int_{\beta_\ringu\R} f(t)g_0(u)d\nu(t,u) = \int f\:g_0\:\nu.
\end{split}
\end{equation}
In the remainder of the paper, we will sometimes use the above right hand side compact notation whenever the variables and domains of integration are clear from the context.

Measures $\nu \in {\mathscr P}(\tspace\times\cuspace)$ obtained as limits of sequences $(u_k)_{k\in\N}\subseteq {\mathscr L}^p(\tspace;\uspace)$ in the sense of \eqref{eq: convergence DPM} have been called {DiPerna-Majda measures}. 
They will be discussed in more detail in Section \ref{sec:DiPerna Majda measures}. It turns out that every DiPerna-Majda measure $\nu \in {\mathscr P}(\tspace\times\cuspace)$ can be disintegrated into a measure $\tau$ on time and an $L^p$-Young measure $\omega$ on $\cuspace$, i.e. $d\nu(t,u)=d\omega(du|t)d\tau(t)$ for some $\tau \in {\mathscr P}(\tspace)$ and $\omega(.|t) \in {\mathscr P}(\cuspace)$. 

A relaxed version of \eqref{eq:classical ocp} taking into account both oscillation and concentration effects can hence be stated as
\begin{equation}
\begin{split}
\inf_{\nu}\int & L_0(t,y(t),u)\,d\nu(t,u)  \\
\st & \int F_0(t,y(t),u)d\nu(t,u) = y_1-y_0, \\
	& \nu \in {\mathscr P}(\tspace\times\cuspace)
\end{split}
\end{equation}
where
\begin{equation}\label{LF}
L_0(t,y,u):=\frac{L(t,y,u)}{1+|u|^p}, \quad
F_0(t,y,u):=\frac{F(t,y,u)}{1+|u|^p}.
\end{equation}
In \cite{Claeys2017}, the moment-sums-of-squares hierarchy is adapted to compute numerically DiPerna-Majda measures and solve optimal control problem featuring oscillations and concentrations.
However, the approach is valid under a certain number of technical assumptions on the data $L$ and $F$, see \cite[Assumption 1, Section 2.2]{Claeys2017}. These assumptions are enforced to prevent the simultaneous presence of concentration and discontinuity.

\subsection{Oscillations, Concentrations and Discontinuities}\label{sec: Oscillation, Concentration, and Discontinuity}
As mentioned in the introduction, the integrals in \eqref{eq:classical ocp} might not be well defined, as concentration effects of the control are likely to cause discontinuities in the trajectory occurring at the same time. In view of the previous examples we propose to generalize the DiPerna-Majda measures, which themselves are a generalization of Young measures, even further and now also relax the trajectory to a measure valued function depending on time and control. In the sequel we describe accordingly the set of {anisotropic parametrized measures}. Then we provide a linear formulation of optimal control problem \eqref{eq:classical ocp} that can cope with oscillations, concentrations and discontinuities in a unified fashion.

\section{Anisotropic Parametrized Measures}\label{sec:extDiPeMa}
In the following we describe the generalized DiPerna-Majda measures. For this it will be instructive to review first the classical DiPerna-Majda measures.

\subsection{DiPerna-Majda measures}\label{sec:DiPerna Majda measures} 
Let $\ringu$ be a complete\footnote{A ring of functions is complete if it contains all constant functions, it separates points from closed subsets and it is closed with respect to the supremum norm.} and separable subring of continuous bounded functions from $\uspace$ to $\R$. It is known \cite[Sect.~3.12.22]{Engelking1989} that there is a one-to-one correspondence between such rings and metrizable compactifications of $\uspace$. By a compactification we mean a compact set, denoted by $\cuspace$, into which $\uspace$ is embedded homeomorphically and densely. For simplicity, we will not distinguish between $\uspace$ and its image in $\cuspace$. Similarly, we will not distinguish between elements of $\ringu$ and their unique continuous extensions defined on $\cuspace$.

DiPerna and Majda \cite{DiPerna1987}, see also \cite{Roubicek1997}, have shown that every bounded sequence $(u_k)_{k\in\N}$ in ${\mathscr L}^p(\tspace;\R^m)$ with $1\le p<\infty$ has a subsequence (denoted by the same indices) such that there exists a probability measure $\tm\in {\mathscr P}(\tspace)$ and an $L^p$-Young measure $\um(.|t)\in {\mathscr P}(\cuspace)$ satisfying for all $f\in {\mathscr C}(\tspace)$ and $g_0\in\ringu$:
\begin{align}\label{eq: convergence DPM disintegrated}
\begin{aligned}
 &\lim_{k\to\infty}\int_0^1 f(t) g_0(u_k(t))(1+\abs{u_k(t)}^p)dt\\ &\qquad = \int_0^1 \int_\cuspace f(t) g_0(u) \um(du|t)\tm(dt) \\
 &\qquad = \int_0^1 \int_\cuspace f(t) g_0(u) d\nu(t,u) = \int f\:g_0\:\nu, 
\end{aligned}  
\end{align}
compare with \eqref{eq: convergence DPM}. The limit measure $d\nu(t,u):=\um(du|t)\tm(dt)$ of such a sequence, or sometimes the pair $(\tm,\um)$, is called a DiPerna-Majda measure.

\subsection{Generalization}
The drawback of DiPerna-Majda measures is that $g$ in \eqref{eq: convergence DPM disintegrated} must be a continuous function. This does not fit to our aim to study interactions of discontinuities and concentrations. 
To go further the simplistic illustration of the introduction, let us consider the following example.

\begin{example}\label{ex1}
Consider  a sequence $(y_k)_{k\in\N}\subset {\mathscr W}^{1,1}([0,1])$ such that $\lim_{k\to\infty}y_k=y$ in ${\mathscr L}^q([0,1])$ for every $1\le q<+\infty$.  We are interested in the integral
$$
\lim_{k\to\infty} \int_{0}^1 g(u_k(t))h(y_k(t))dt
$$
for continuous functions $g$ and $h$  such that  $|g(u)|\le C(1+|u|)$ with some constant $C>0$, and where $u_k:=\dot{y}_k \in {\mathscr L}^1([0,1])$ is the weak derivative of $y_k$. If $g$ is the identity then the calculation is easy, namely
the limit equals $\liminf_{k\to\infty} H(y_k(1))-H(y_k(0))$ where $H$ is the primitive of $h$. In the case of a more general function $g$, the situation is more involved. For example for $k\geq 2$ let 
$$
u_k(t):=
\begin{cases}
0 &\text{ if $0\le t\le \frac{1}{2}$},\\
k &\text{ if $\frac{1}{2}\le t\le \frac{1}{2}+\frac{1}{k}$},\\
0 &\text{ if $\frac{1}{2}+\frac{1}{k}\le t\le 1$}
\end{cases}
$$
whose primitive is
$$
y_k(t):=
\begin{cases}
0 &\text{ if $0\le t\le \frac{1}{2}$},\\
k(t-\frac{1}{2}) &\text{ if $\frac{1}{2}\le t\le \frac{1}{2}+\frac{1}{k}$},\\
1 &\text{ if $\frac{1}{2}+\frac{1}{k}\le t\le 1$}
\end{cases}
$$
see Figure \ref{fig:ex1}.
\begin{figure}[h]
\begin{center}
\includegraphics[scale=.5]{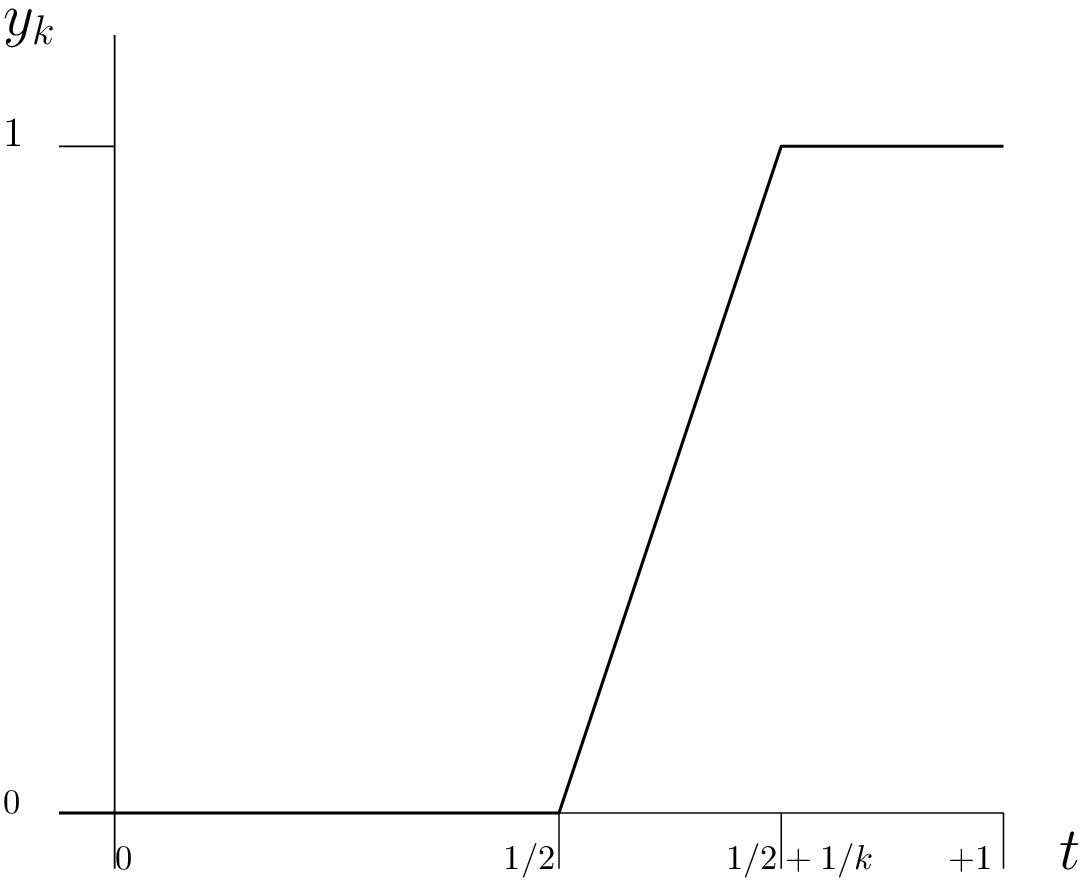}
\includegraphics[scale=.5]{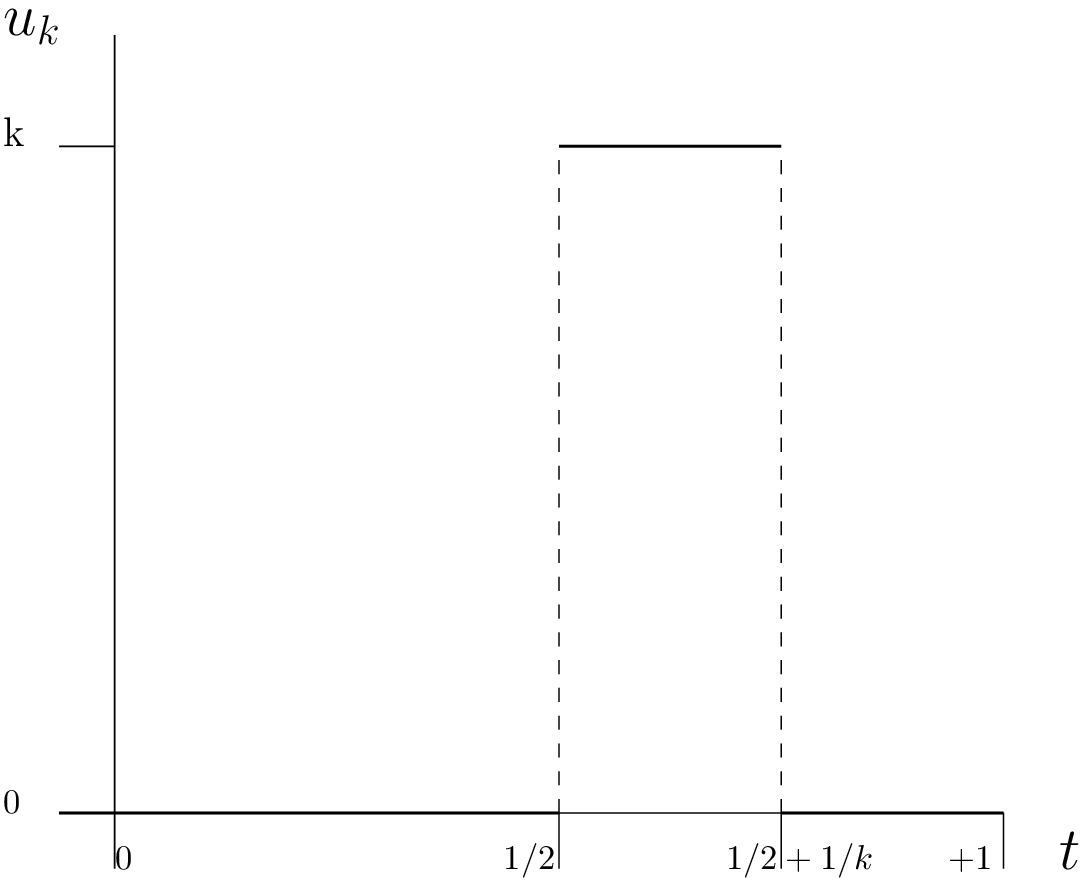}
\caption{Sequences $(y_k,u_k)_{k\in\N}$ from Example~\ref{ex1}.\label{fig:ex1}}
\end{center}
\end{figure}
It is easy to see that
\[
\begin{array}{l}
\lim_{k\to\infty} \int_0^1 g(u_k(t))h(y_k(t))dt
= \int_0^{\frac{1}{2}} g(0)h(0)dt + \\ \quad
\lim_{k\to\infty} \int_{\frac{1}{2}}^{\frac{1}{2}+\frac{1}{k}} g(k)h(k(t-\frac{1}{2}))dt + \\ \quad
\lim_{k\to\infty} \int_{\frac{1}{2}+\frac{1}{k}}^1 g(0)h(1)dt
= \frac{1}{2}g(0)(h(0)+h(1))+ \\ \quad
\lim_{k\to\infty} \int_{\frac{1}{2}}^{\frac{1}{2}+\frac{1}{k}} \frac{\dot{H}(k(t-\frac{1}{2}))}{k}g(k)dt
=\frac{1}{2}g_0(0)(h(0)+h(1))+ \\ \quad 
(H(1)-H(0))\lim_{k\to\infty}\frac{g(k)}{k}.
\end{array}
\]
The sequence $(u_k)_{k\in\N}$ concentrates at $\frac{1}{2}$ which is exactly the point of discontinuity of the pointwise limit of $(y_k)_{k\in\N}$.
Also notice that  $u_k$ converges weakly to $\delta_{\frac{1}{2}}$ in ${\mathscr P}([0,1])$ when $k\to\infty$. \new{The factor $H(1)-H(0)$ in the previous equation suggests that  we should refine the definition of the pointwise limit of $(y_k)_{k\in\N}$ at $\frac{1}{2}$ by enforcing that it is the  Lebesgue measure supported  on the interval of the jump. We will make this rigourous in the following. The other term in the factor also shows that the limit of $g(k)/k$ should exist when $k$ tends to infinity.}
\end{example}

To cope with the simultaneous presence of oscillations, concentrations and discontinuities, a new tool was recently introduced in \cite{Kalamajska2017}, namely anisotropic parametrized measures generated by pairs  $(y_k,u_k)_{k\in\N}$ where $u_k$ is the control and $y_k$ the corresponding state trajectory. Let us describe now what we need in our optimal control context. First, let us make the following observation:

\begin{lemma}\label{ybounded}
Any admissible trajectory of optimal control problem (\ref{eq:classical ocp}) is such that $y \in {\mathscr L}^\infty([0,1];Y)$ for some compact set $Y \subset \R^n$, e.g. a ball of sufficiently large radius.
\end{lemma}

{\bf Proof:}
The function $t \mapsto y(t)$ is the integral of a Lesbesgue integrable function, and on a bounded time interval, it is bounded.$\Box$

Then, the following result is a special case of \cite[Theorem~2]{Kalamajska2017}:

\begin{theorem}\label{thm1}
Let  $1\le p< +\infty$. Let  $(u_k)_{k\in\N}$ be a bounded sequence in  ${\mathscr L}^p(\tspace;\R^m)$ and $(y_k)_{k\in\N}$ a bounded sequence in ${\mathscr W}^{1,1}(\tspace;\R^n)$. Then  there is a (non-relabeled) subsequence $(u_k,y_k)_{k\in\N}$,  a measure $\tm \in {\mathscr P}(\tspace)$, a measure $\um(.|t) \in {\mathscr P}(\cuspace)$ parametrized in $t\in\tspace$ and a measure $\upsilon(.|t,u) \in {\mathscr P}(Y)$ parametrized in $t\in\tspace$ and $u\in\cuspace$ such that for every $f \in {\mathscr C}(\tspace)$, $g_0\in\ringu$, $h_0\in {\mathscr C}(Y)$, it holds
\begin{align}\label{generate}
\begin{aligned}
&\lim_{k\to\infty}\int_0^1 f(t)g_0(u_k(t))(1+|u_k(t)|^p)h_0(y_k(t))dt\\
&=\int_0^1 \int_{\cuspace}\int_{Y}f(t)g_0(u)h_0(y)\upsilon(dy|t,u)\omega(du|t)\tau(dt)\\
&=\int_0^1 \int_{\cuspace}\int_{Y} f(t)g_0(u)h_0(y)d\mu(t,y,u) \\
&= \int f\:g_0\:h_0\:\mu.\\
\end{aligned}
\end{align}
The measure $d\mu(t,u,y):=\upsilon(dy|t,u)\omega(du|t)\tau(dt)$, or sometimes the triplet $(\tau,\omega,\upsilon)$, is called an anisotropic parametrized measure. Moreover, the \new{DiPerna-Majda} measure $(\tau,\omega)$ is generated by 
$(u_k)_{k\in\N}$.
\end{theorem}

\begin{example}\label{ex1b}
Let us revisit Example \ref{ex1} and the calculations of the integrals there. Let $f \in {\mathscr C}([0,1])$, let $h \in {\mathscr C}(\R)$ be bounded with primitive denoted by $H$, and let $g:=(1+|.|)g_0$ where $g_0 \in \ringu$ corresponding to the two-point (or sphere) compactification $\cuspace = \R \cup \{\pm\infty\}$, i.e. such that $\lim_{u\to\pm\infty}g_0(u)=:g_0(\pm\infty)\in\R$. Then it holds
\[
\begin{array}{l}
\lim_{k\to\infty} \int_0^1 f(t)g(u_k(t))h(y_k(t))dt
= \int_0^{\frac{1}{2}} f(t)g(0)h(0)dt + \\ \quad
\lim_{k\to\infty} \int_{\frac{1}{2}}^{\frac{1}{2}+\frac{1}{k}} f(t)g(k)h(k(t-\frac{1}{2})))dt + \\ \quad
\lim_{k\to\infty} \int_{\frac{1}{2}+\frac{1}{k}}^1 f(t)g(0)h(1)dt
= \int_0^{\frac{1}{2}} f(t)g(0)h(0)dt  + \\ \quad
 \int_{\frac{1}{2}}^1 f(t)g(0)h(1)dt + \lim_{k\to\infty} \int_{\frac{1}{2}}^{\frac{1}{2}+\frac{1}{k}} f(t)g(k)\frac{\dot{H}(k(t-\frac{1}{2}))}{k}dt \\ \quad
= \int_0^{\frac{1}{2}} f(t)g(0)h(0)dt  + \int_{\frac{1}{2}}^1 f(t)g(0)h(1)dt + \\ \quad
 \lim_{k\to\infty} \int_{\frac{1}{2}}^{\frac{1}{2}+\frac{1}{k}} f(t)g_0(k)\dot{H}(k(t-\frac{1}{2}))\frac{1+k}{k}dt \\ \quad
= \int_0^{\frac{1}{2}} f(t)g(0)h(0)dt  + \int_{\frac{1}{2}}^1 f(t)g(0)h(1)dt + \\ \quad
f(\frac{1}{2})g_0(+\infty)(H(1)-H(0))  \\ \quad
=\int_0^1 \int_\cuspace \int_Y f(t)g_0(u)h(y)\upsilon(dy|t,u)\omega(du|t)\tau(dt)
\end{array}
\]
where
\[
\tau(dt) = \lambda_{[0,1]} + 2\delta_{\frac{1}{2}}
\]
and
\[
\omega(du|t) = \begin{cases}
\delta_{+\infty} & \text{ if $t=\frac{1}{2}$},\\
\delta_0 & \text{ otherwise}
\end{cases}
\]
and
$$
\upsilon(dy|t,u) =\begin{cases}
\delta_0 &\text{ if $t\in[0,\frac{1}{2})$ },\\
\lambda_{[0,1]} & \text{ if $t=\frac{1}{2}$},\\
\delta_1 &\text{ if $t\in(\frac{1}{2},1]$ }\\
\end{cases}
$$
where $\lambda_X$ denotes the Lebesgue measure on $X$, and $Y = [0,1]$.
\end{example}
\begin{example}\label{ex2}
Let us revisit the slightly more complicated \cite[Example 3]{Kalamajska2017}, appropriately scaled on $[0,1]$. The trajectory sequence is 
$$
y_k(t):=\begin{cases}
0 &\text{ if $0\le t\le \frac{1}{2}-\frac{1}{k}$},\\
k(t-\frac{1}{2}+\frac{1}{k}) &\text{ if $\frac{1}{2}-\frac{1}{k} \le t \le \frac{1}{2}$},\\
-2k(t-\frac{1}{2}-\frac{1}{2k}) &\text{ if $\frac{1}{2} \le t\le \frac{1}{2}+\frac{1}{k}$},\\
-1 & \text{ if $\frac{1}{2}+\frac{1}{k}\le t\le 1$}
\end{cases}
$$
and its weak derivative $u_k:=\dot{y}_k$ is
$$
u_k(t):=\begin{cases}
0 &\text{ if $0\le t\le \frac{1}{2}-\frac{1}{k}$},\\
k &\text{ if $\frac{1}{2}-\frac{1}{k} \le t \le \frac{1}{2}$},\\
-2k &\text{ if $\frac{1}{2} \le t\le \frac{1}{2}+\frac{1}{k}$},\\
0 & \text{ if $\frac{1}{2}+\frac{1}{k}\le t\le 1$}
\end{cases}
$$
see Figure \ref{fig:ex2}. 
\begin{figure}[h]
\begin{center}
\includegraphics[scale=.5]{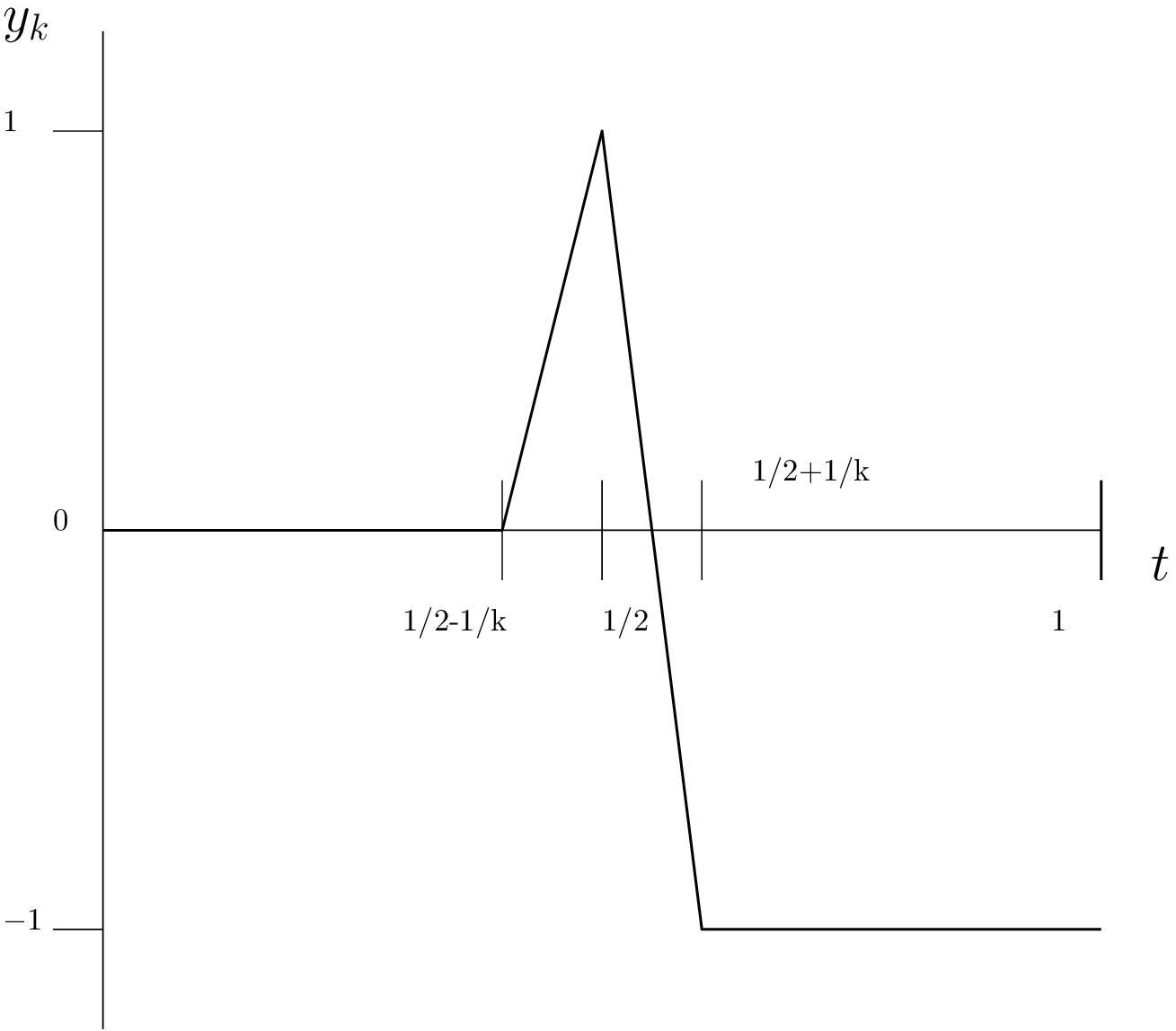}
\includegraphics[scale=.5]{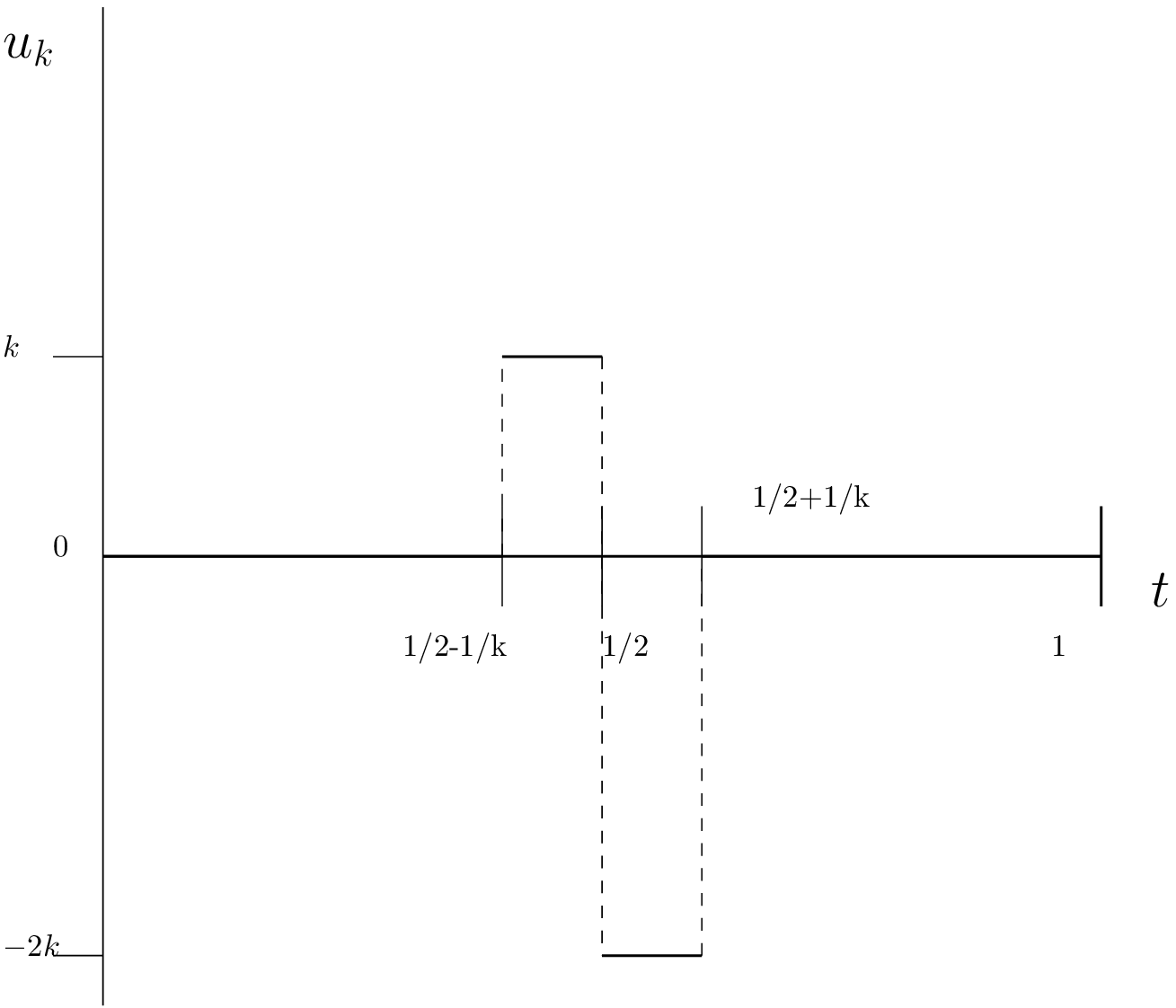}
\caption{Sequences $(y_k,u_k)_{k\in\N}$ from Example~\ref{ex2}.}
\end{center}
\label{fig:ex2}
\end{figure}
 Let $f \in {\mathscr C}([0,1])$, let $h \in {\mathscr C}(\R)$ be bounded with primitive denoted by $H$, and let $g=(1+|.|)g_0$ where $g_0 \in \ringu$ corresponding to the two-point (or sphere) compactification $\cuspace = \R \cup \{\pm\infty\}$, i.e. such that $\lim_{u\to\pm\infty}g_0(u)=:g_0(\pm\infty)\in\R$. Then it holds
\[
\begin{array}{l}
\lim_{k\to\infty} \int_0^1 f(t)g(u_k(t))h(y_k(t))dt \\ \quad
= \lim_{k\to\infty} \int_0^{\frac{1}{2}-\frac{1}{k}} f(t)g(0)h(0)dt + \\ \quad 
\lim_{k\to\infty} \int_{\frac{1}{2}-\frac{1}{k}}^{\frac{1}{2}} f(t)g(k)h(k(t-\frac{1}{2}+\frac{1}{k}))dt + \\ \quad
 \lim_{k\to\infty} \int_{\frac{1}{2}}^{\frac{1}{2}+\frac{1}{k}} f(t)g(-2k)h(-2k(t-\frac{1}{2}-\frac{1}{2k}))dt + \\ \quad
\lim_{k\to\infty} \int_{\frac{1}{2}+\frac{1}{k}}^1 f(t)g(0)h(-1)dt
 = \int_0^{\frac{1}{2}} f(t)g(0)h(0)dt + \\ \quad
\int_{\frac{1}{2}}^1 f(t)g(0)h(-1)dt +\\ \quad
 \lim_{k\to\infty}  \int_{\frac{1}{2}-\frac{1}{k}}^{\frac{1}{2}} f(t)g_0(k)\dot{H}(k(t-\frac{1}{2}+\frac{1}{k}))\frac{1+k}{k}dt + \\ \quad
 \lim_{k\to\infty} \int_{\frac{1}{2}}^{\frac{1}{2}+\frac{1}{k}} f(t)g_0(-2k)\dot{H}(-2k(t-\frac{1}{2}-\frac{1}{2k}))\frac{1+2k}{-2k}dt \\ \quad
= \int_0^{\frac{1}{2}} f(t)g(0)h(0)dt + \int_{\frac{1}{2}}^1 f(t)g(0)h(-1)dt + \\ \quad
f(\frac{1}{2})g_0(+\infty)(H(1)-H(0)) + f(\frac{1}{2})g_0(-\infty)(H(1)-H(-1))  \\ \quad
=\int_0^1  \int_\cuspace \int_Y f(t)g_0(u)h(y)\upsilon(dy|t,u)\omega(du|t)\tau(dt)
\end{array}
\]
where
\[
\tau(dt) = \lambda_{[0,1]} + 3\delta_{\frac{1}{2}}
\]
and
\[
\omega(du|t) = \begin{cases}
\frac{1}{2}\delta_{+\infty}+\frac{1}{2}\delta_{-\infty} & \text{ if $t=\frac{1}{2}$},\\
\delta_0 & \text{ otherwise}
\end{cases}
\]
and
$$
\upsilon(dy|t,u) =\begin{cases}
\delta_0 &\text{ if $t\in[0,\frac{1}{2})$ },\\
\lambda_{[0,1]} & \text{ if $t=\frac{1}{2}$ and $u=+\infty$},\\
\frac{1}{2}\lambda_{[-1,1]} & \text{ if $t=\frac{1}{2}$ and $u=-\infty$},\\
\delta_{-1} &\text{ if $t\in(\frac{1}{2},1]$ }\\
\end{cases}
$$
where $\lambda_X$ denotes the Lebesgue measure on $X$, and $Y=[-1,1]$.
\end{example}

\section{Relaxed Optimal Control with Oscillations, Concentrations and Discontinuities}\label{sec:OCD}

To the classical optimal control problem (\ref{eq:classical ocp}) we associate the relaxed optimal control problem
\begin{equation}\label{rocp}
\begin{split}
v^*_R := \inf_\mu \int & L_0 \: \mu \\
\st & \int F_0 \: \mu = y_T - y_0,\\
	&  \mu \in {\mathscr P}(\tspace\times\cuspace\times Y)
\end{split}
\end{equation}
which is {\sl linear} in the unknown  measure $\mu$. In contrast, classical problem (\ref{eq:classical ocp}) is non-linear in the unknown trajectory $y$ and control $u$.

Since optimal control problem (\ref{rocp}) is a relaxation of the optimal control problem (\ref{eq:classical ocp}), it may happen that the infimum in (\ref{rocp}) is strictly less than the infimum in (\ref{eq:classical ocp}), i.e. $v^*_R < v^*$. Formulating necessary and sufficient conditions on the problem data $F$ and $L$ such that $v^*_R = v^*$, i.e. there is no relaxation gap is an open problem.  However, if we know that the \new{probability} measure in problem (\ref{rocp}) is generated by limits of functions, then there is no relaxation gap. Let us explain this now.

\begin{assumption}[Regularity of the data]\label{regularity}
Let $L$ and $F$ be such that in (\ref{LF}) it holds
\begin{align}\label{L0}
L_0\in {\mathscr C}([0,1]\times\cuspace\times Y)
\end{align}
and
\begin{align}\label{F0}
F_0\in {\mathscr C}([0,1]\times\cuspace\times Y;\R^n).
\end{align}
Moreover, there is a constant $c_L>0$ such that
\begin{align}\label{L}
L(t,u,y)\ge c_L|u|^p
\end{align}
for all $t$, $u$, $y$ and there is a constant $c_F>0$ such that
\begin{align}\label{lipschitz}
|F(t,u,y_1)-F(t,u,y_2)|\le c_F(|u|^p+1)|y_1-y_2|
\end{align}
for all $t$, $u$, $y_1$, $y_2$.
\end{assumption}

The following result follows from \new{classical existence and uniqueness results for differential equations, see e.g. \cite[Theorem 3.1]{Bressan2007}:}
\begin{lemma}\label{cara}
Assume that $p\ge 1$, $u\in {\mathscr L}^p([0,1];\R^m)$  and  $y_0\in\R^n$ are given. Let further $F: [0,1]\times\R^m\times\R^n\to\R^n$ satisfy \eqref{F0} and \eqref{lipschitz}. Then  
\begin{align}\label{equation}
dy(t)=F(t,u(t),y(t))dt \ , \ y(0)=y_0
\end{align}
has a unique solution $y \in {\mathscr L}^\infty(\tspace;Y)$ with values in a compact subset $Y$ of $\yspace$.
\end{lemma}

Assume that there is a bounded sequence  $(u_k)_{k\in\N}\subset {\mathscr L}^p$  and that $\{y_k\}_{k\in\N} \subset {\mathscr W}^{1,1}$ is a sequence of corresponding solutions obtained in Lemma~\ref{cara}. Then $\{y_k\}$ is uniformly bounded in ${\mathscr W}^{1,1}$. Indeed, due to \eqref{F0} 
we  see that $\frac{d|y_k(t)|}{dt} \le \left|\frac{dy_k(t)}{dt}\right| =|F(t,u_k(t),y_k(t))|\le c_F(1+|u_k(t)|^p+|y_k(t)|)$. Then the Gronwall inequality \cite[Appendix B.2.j]{Evans98} implies that $\sup_{k\in\N}\|y_k\|_{W^{1,1}} < \infty$
and since $y_k$ is the integral of an integrable function on a bounded time interval, it holds that $\new{y_k} \in {\mathscr L}^\infty(\tspace;Y)$ for $Y \subset \yspace$ a ball of radius $\sup_{k\in\N}\|y_k\|_{L^\infty} < \infty$.
The limit of the right-hand side of \eqref{equation} can then be expressed in terms of an anisotropic parametrized measure $\mu$:
\begin{equation}\label{limitF}
\lim_{k\to\infty} F(t,u_k(t),y_k(t))dt =  \int_{\beta_\mathcal{U}\R^{m}}\int_{Y}F_0(t,u,y) d\mu(t,u,y).
\end{equation}

\new{
As explained in \cite[Theorem 7]{Kalamajska2017}, the integral (\ref{generate}) in the definition of the anisotropic parametrized measure can be decomposed as follows
\begin{equation}\label{decomp}
\begin{array}{l}
\displaystyle
\int_0^1 \int_{\cuspace}\int_{Y} f(t)g_0(u)h_0(y)d\mu(t,y,u) =\\
\quad \displaystyle \int_0^1 \int_{\R^m} f(t)g_0(u)(1+|u|^p)h_0(y(t))\tilde{\omega}(du|t)dt + \\
\quad \displaystyle  \int_0^1 \int_{\cuspace\setminus\R^m}\int_{Y}f(t)g_0(u)h_0(y)\upsilon(dy|t,u)\omega(du|t)\tau(dt)
\end{array}
\end{equation}
where $\tilde{\omega}$ is a classical Young measure on $\R^m$.
Using the decomposition (\ref{decomp}), instead of (\ref{equation}) we get the following differential equation
\begin{equation}\label{eqn-meas}
\begin{array}{l}
dy(t) = F(t,u,y(t))\tilde{\omega}(du|t)dt + \\ \qquad\displaystyle \int_{\beta_\mathcal{U}\R^{m}\setminus\R^{m}} \int_Y F_0(t,u,y)\upsilon(dy|t,u)\omega(du|t)\tau(dt).
\end{array}
\end{equation}

}

\begin{lemma}
Given an anisotropic parametrized measure $\mu$ and an initial condition $y_0$, the solution $y$ to \eqref{eqn-meas} is unique.
\end{lemma}

{\bf Proof:}
Assume that it is not the case, i.e., that there are two solutions $y_1, y_2\in {\mathscr L}^\infty(\tspace;Y)$. Desintegrating $d\mu(t,y,u)=\upsilon(dy|t,u)\omega(du|t)\tau(dt)$, we get the following relationship for the difference $y_d:=y_1-y_2$ because of \eqref{lipschitz}, it holds
$|\dot y_d|\le \int_{\R^m}|F(t,u,y_1(t))-F(t,u,y_2(t))|\omega_t(du) \le \int_{\R^m}c_F(|u|^p+1)\omega_t(du)|y_d(t)|$.
The right hand side belongs to ${\mathscr L}^1(\tspace)$, therefore the measure $dy_d(t)$ is absolutely continuous with respect to the uniform measure $dt$. As $y_d(0)=0$  we have $y_d(t)=0$ for all $t \in \tspace$, by the Gronwall inequality \cite[Appendix B.2.j]{Evans98}.$\Box$

In relaxed optimal control problem (\ref{rocp}) we use an integral formulation of (\ref{eqn-meas}) incorporating the initial and terminal conditions: $\int_0^1 \int_{\cuspace} \int_{Y} F_0(t,u,y)d\mu(t,u,y) = \int F_0\:\mu = y_1 - y_0$.
For each anisotropic parametrized measure $\mu$, we can therefore associate a sequence of trajectories  $\{y_k\} \subset {\mathscr W}^{1,1}$ and controls $(u_k) \subset {\mathscr L}^p$ satisfying the differential equation (\ref{equation}) and such that (\ref{limitF}) holds. Conversely, the limit of each such sequence of trajectories and controls
can be modeled by an anisotropic parametrized measure. 
The following result of absence of relaxation gap then follows immediately from
the construction of problem (\ref{rocp}).

\begin{proposition}[No relaxation gap]\label{norelaxationgap}
Let Assumption \ref{regularity} hold \new{and let $\mu$ solve problem (\ref{rocp}).
If there is an admissible sequence $(u_k, y_k)_{k\in\N}$ such that (\ref{generate}) holds}
then $v^*_R=v^*$.
\end{proposition}

\section{Relaxed Optimal Control with Occupation Measures}\label{sec:relatcontrol} 

In the previous section, we proposed a linear reformulation of non-linear optimal control, thanks to the introduction of anisotropic parametrized measures. In the current section, we describe another linear reformulation proposed in \cite{Lasserre2008} and relying on the notion of occupation measure. The relation between this linear reformulation and the classical Majda-DiPerna measures was investigated in \cite{Claeys2017}, with the help of a graph completion argument. In the sequel we show that the generalized Majda-DiPerna measures also fit naturally this framework.

Let $v \in {\mathscr C}^1(\tspace\times Y)$. For any admissible trajectory $y$ and control $u$ solving the differential equation (\ref{equation}), it holds
$\int_0^1 dv(t,y(t)) = v(1,y(1))-v(0,y(0)) = \int_0^1 \left(\frac{\partial v}{\partial t}(t,y(t)) + \frac{\partial v}{\partial y}(t,y(t)) \cdot \dot{y}(t)\right) dt$.
Optimal control problem (\ref{eq:classical ocp}) can then be rewritten as
\begin{equation}\label{eq:weak-ocp}
\begin{split}
v^* = \inf_u & \int_0^1 L(t,u(t),y(t)) dt\\
\st & \int_0^1 \left(\frac{\partial v}{\partial t} + \frac{\partial v}{\partial y} \cdot F\right)(t,u(t),y(t)) dt \\
& = v(1,y_1)-v(0,y_0), \: \forall v \in {\mathscr C}^1(\tspace\times\R^n)\\
	& y\in {\mathscr W}^{1,1}(\tspace;\R^n), \quad u\in {\mathscr L}^p(\tspace;\R^m).
\end{split}
\end{equation}

\begin{definition}[Occupation measure]\label{def:occup-meas}
Given a control $u$ and a trajectory $y$ solving the differential equation (\ref{equation}), we define the occupation measure $\mu_{u,y} \in y\in {\mathscr P}(\tspace\times\R^n\times\R^m)$ by
\[
d\mu_{u,y}(t,u,y):=\delta_{y(t)}(dy)\delta_{u(t)}(du)dt.
\]
\end{definition}
Geometrically $\mu_{u,y}(A\times B\times C)$ is the time spent by the trajectory $(t,u(t),y(t))$ in any Borel subset $A\times B\times C$ of $\tspace\times\R^m\times Y$. Analytically, integration with respect to $\mu_{u,y}$ is the same as integration along $(u(t),y(t))$ with respect to time. In particular
$\int_0^1 L(t,u(t),y(t))dt = \int_0^1 \int_{\R^m} \int_{\R^n} L(t,u,y)d\mu_{u,y}(t,u,y) = \int L \:\mu_{u,y}$
and for all test functions $v \in {\mathscr C}^1(\tspace\times Y)$, it holds that
$\int_0^1  \left(\frac{\partial v}{\partial t} + \frac{\partial v}{\partial y} \cdot F\right)(t,u(t),y(t)) dt
$\\$= \int_0^1 \int_{\R^m} \int_{Y} \left(\frac{\partial v}{\partial t} + \frac{\partial v}{\partial y} \cdot F\right)(t,u,y) d\mu_{u,y}(t,u,y)
 = \int  \left(\frac{\partial v}{\partial t}  + \frac{\partial v}{\partial y}  \cdot F \right) \mu_{u,y}$.
Using the same arguments as in \cite[Proposition 4]{Claeys2017}, we can reformulate  optimal control problem \eqref{eq:weak-ocp} as a linear problem on measures, leading to the following relaxed formulation:
\begin{equation}\label{eq:relat-ocp}
\begin{split}
v^*_M := \inf_\mu & \int L_0\:\mu \\
\st & \int \left(\frac{\partial v}{\partial t}(1+|u|^p)^{-1} + \frac{\partial v}{\partial y} \cdot F_0\right) \mu \\
& = v(1,y_1)-v(0,y_0) \quad \forall v \in {\mathscr C}^1(\tspace\times Y),\\
& \mu \in {\mathscr P}(\tspace\times\cuspace\times Y).
\end{split}
\end{equation}
Note that $\mu$ in the above problem is not necessarily an occupation measure in the sense of Definition \ref{def:occup-meas}, but a general probability measure in ${\mathscr P}(\tspace\times\cuspace\times Y)$. For this reason, the infimum in relaxed problem (\ref{eq:relat-ocp}) can be strictly less than the infimum in classical problem (\ref{eq:classical ocp}), i.e. $v^*_M < v^*$.
	
\begin{proposition}[No relaxation gap]
It holds $v^*_R \leq v^*_M \leq v^*$ and hence if there is no relaxation gap in relaxed problem (\ref{rocp}) then there is no relaxation gap in relaxed problem (\ref{eq:relat-ocp}).
\end{proposition}

{\bf Proof:}
Just observe that problem (\ref{rocp}) corresponds to the particular choice of test functions $v(t,y):=y_k$, $k=1,\ldots,n$ in problem (\ref{eq:relat-ocp}). Hence the infimum in (\ref{rocp}) is smaller than the infimum in (\ref{eq:relat-ocp}), which is in turn smaller than the infimum in (\ref{eq:classical ocp}), i.e. $v^*_R \leq v^*_M$. Now if $v^*_R=v^*$ then obviously $v^*_M=v^*$. $\Box$

\section{\new{The Lasserre hierarchy}}\label{sec:approach}

Once we get to the measure linear problem \eqref{eq:relat-ocp}, we follow the same strategy as in \cite[Section 4]{Claeys2017}:
\begin{enumerate}
 \item compactify the control space by using a change of variables and homogenization;
  \item since all the data are polynomial, construct an equivalent moment linear problem where the unknown are moments of the occupation measure supported on a compact semialgebraic set;  
 \item use the moment-sums-of-squares hierarchy as in \cite{Lasserre2008} to obtain a sequence of approximate moments at the price of solving numerically semidefinite programming problems; 
 \item from the approximate moments, construct an approximate solution to the optimal control problem. 
\end{enumerate}

\section{\new{Illustrative example}}\label{sec:example}

Let us illustrate this strategy on our introductory example (\ref{ex}). The trajectory $y$ should move the state from zero at initial time to one at final time, yet for the non-negative integrand to be as small as possible, the control $u$ should be zero all the time, except maybe at time zero. 
If $\varepsilon=1$ this problem has a trivial optimal solution $u(t)=0$. For $\varepsilon=0$ as explained already we can solve the problem by integration by parts because $\dot{y}(t) = u(t)$. The integration trick cannot be carried out in the case of $\varepsilon\in(0,1)$. 

We use the relaxation proposed in Section \ref{sec:relatcontrol} to formulate problem \eqref{ex} as a measure LP:
\begin{equation}\label{eq:mlp1}
\begin{split}
\inf_{\mu} & \int \frac{(t+y)u}{1+u}\mu \\
\st & \int  \left(\frac{\partial v}{\partial t} \frac{1}{1+u} + \frac{\partial v}{\partial y} \frac{\sqrt{\varepsilon^2+u^2}}{1+u}\right) \mu  \\
& = v(1,1)  - v(0,0), \text{ for all } v\in{\mathscr C}^1([0,1]^2)\\
	& \mu\in{\mathscr P}([0,1]\times\beta\R_+\times[0,1]).
\end{split}
\end{equation}
Note that we can omit the absolute value in the denominator, as $u$ is constrained to be non-negative. 

\new{Since in problem (\ref{ex}) the growth of the Lagrangian and the dynamics is at most linear,}
we expect the control to concentrate. Therefore let $u(t):=\tfrac{r(t)}{1-r(t)}$ with $r(t)\in[0,1]$. Then the dynamic of $y$ reads $\dot{y}(t) = \sqrt{\left(\tfrac{r(t)}{1-r(t)}\right)^2 + \varepsilon^2} = \frac{\sqrt{r(t)^2+\varepsilon^2(1-r(t))^2}}{1-r(t)}$.
Introduce the auxiliary variable $w(t)$ such that $w(t)^2 = r(t)^2+\varepsilon^2(1-r(t))^2$. By knowledge of bounds for $\varepsilon$ and $r(t)$ we can conclude that $0\leq w(t)\leq 1$. The linear problem on moments \new{then} reads
\begin{equation}\label{eq:bmlp1}
\begin{split}
\inf_{\gamma} & \int \left(t +y\right)r\, \gamma\\
\st &  \int \left(\frac{\partial v}{\partial t} (1-r) + \frac{\partial v}{\partial y} w \right) \gamma  \\
& =v(1,1)  - v(0,0), \text{ for all } v\in\R[t,y],\\
	& \gamma\in{\mathscr P}([0,1]^3).\\
\end{split}
\end{equation}

With the following GloptiPoly script we could solve the problem numerically for different values of the parameter $\varepsilon$. \new{From these numerical solutions} we could guess the analytic optimal solution.

\begin{figure}[h]
\begin{center}
\includegraphics[width=\columnwidth]{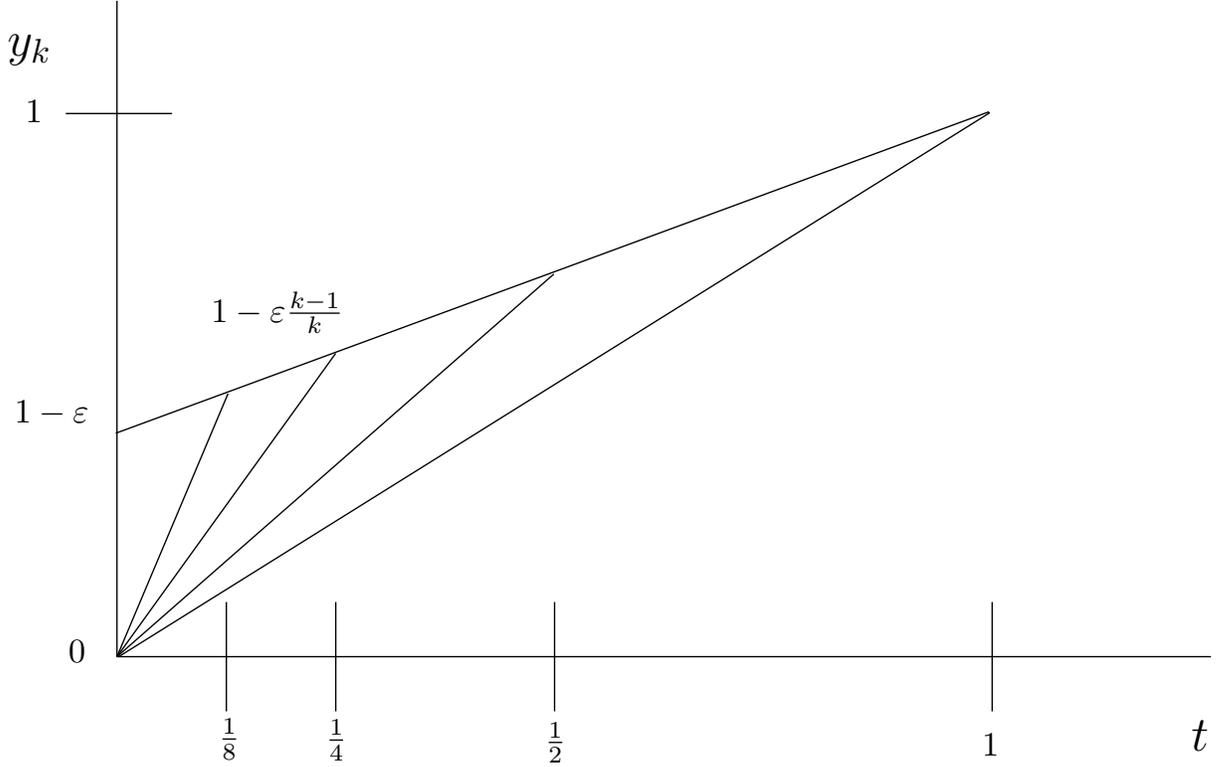}
\caption{Sequence $(y_k)_{k=1,2,4,8}$ from Example~\ref{ex}.\label{fig:ex}}
\end{center}
\end{figure}
The measure $d\mu(t,y,u) = \tau(dt)\omega(du|t)\upsilon(dy|t,u)$ with
\begin{eqnarray}
\tau(dt) = \lambda_{[0,1]} + (1-\varepsilon)\delta_0\\
\omega(du|t)=\left\lbrace\begin{array}{cl}\delta_\infty, & t=0\\ \delta_0,&t>0 \end{array}\right.\\
\upsilon(dy|t,u)=\left\lbrace\begin{array}{cl}\tfrac{1}{1-\varepsilon}\lambda_{[0,1-\varepsilon]}, & t=0\\\delta_{1-\varepsilon+\varepsilon t},&t>0 \end{array}\right.
\end{eqnarray}
\new{solves the relaxation (\ref{rocp}) and hence yields a lower bound  $\tfrac{(1-\varepsilon)^2}{2}$ on the optimum}.
\new{Moreover, this optimum is attained by the sequences}
\[
u_k(t) = \left\lbrace\begin{array}{cl}\sqrt{(k(1-\varepsilon)+\varepsilon)^2-\varepsilon^2}, & t\in[0,\tfrac{1}{k}] \\ 0,&t>\tfrac{1}{k} \end{array}\right.,
\]
and
\[
y_k(t) = \left\lbrace\begin{array}{cl} (k(1-\varepsilon)+\varepsilon)t, & t\in[0,\tfrac{1}{k}] \\ \varepsilon t+ 1 - \varepsilon,&t>\tfrac{1}{k} \end{array}\right.
\]
see Figure \ref{fig:ex}, \new{which proves that we got the optimal solution according to Proposition \ref{norelaxationgap}}.

\begin{figure}
\begin{verbatim}
d = 6; % relaxation order (moments of degree <= 2*d)
e = 0.2; % epsilon

mpol t y r w
gamma = meas([t y r w]); % measure
 
v = mmon([t y],2*d);  % test functions

assign([t y r w],[0 0 0 0]);
v0 = double(v); % initial condition			

assign([t y r w],[1 1 0 0]);	
vT = double(v); % terminal condition			

% moment problem
obj = min((t+y)*r);
spt = [t-t^2>=0, y-y^2>=0, r-r^2>=0, w-w^2>=0, ..
                   w^2==r^2+e^2*(1-r)^2];
dyn = [vT-v0 == mom(diff(v,t)*(1-r) + diff(v,y)*w)];
P = msdp(obj, spt, dyn);

% solve semidefinite relaxation
[stat,bnd] = msol(P);

% display moments of solution
if stat>=0
  double(mom(mmon(t,2*d))) % time
  double(mom(mmon(y,2*d))) % trajectory
  double(mom(mmon(r,2*d))) % normalized control
end
\end{verbatim}
\caption{GloptiPoly script.}
\label{fig:code-compact-time}
\end{figure}

The numerical methods obtained with \new{the GloptiPoly script of Figure \ref{fig:code-compact-time}} and the SeDuMi semidefinite solver for the 6th relaxation (i.e. moments of degree up to 12) are reported in Table \ref{tab: gloptipoly solution}. They match to 4 significant digits with the analytic moments reported in Table \ref{tab: analytic moments}.

\begin{table}
\begin{center}
\begin{tabular}{r|cccc}
$k$   & $\int t^k\mathrm{d}\mu$ & $\int y^k\mathrm{d}\mu$ & $\int r^k\mathrm{d}\mu$ & $\int w^k\mathrm{d}\mu$\\
\hline
 	 0&  	1.8000  &  1.8000  &  1.8000  &  1.8000\\
 	 1&  	0.5000  &  1.2200  &  0.8000  &  1.0000\\
 	 2&    	0.3333  &  0.9840  &  0.8000  &  0.8400\\
 	 3&   	0.2500  &  0.8404  &  0.8000  &  0.8080\\
 	 4&   	0.2000  &  0.7379  &  0.8000  &  0.8016\\
 	 5&    	0.1667  &  0.6586  &  0.8000  &  0.8003\\
 	 6&    	0.1429  &  0.5944  &  0.8000  &  0.8001\\
 	 7&    	0.1250  &  0.5411  &  0.8000  &  0.8000\\
	 8&    	0.1111  &  0.4959  &  0.8000  &  0.8000\\
	 9&    	0.1000  &  0.4571  &  0.8000  &  0.8000\\
	10&    	0.0909  &  0.4233  &  0.8000  &  0.8000\\
	11&    	0.0833  &  0.3938  &  0.8000  &  0.8000\\
	12&    	0.0769  &  0.3677  &  0.8000  &  0.8000
\end{tabular}
\caption{Approximate moments for $\varepsilon=0.2$, computed with Gloptipoly and SeDuMi.}
\label{tab: gloptipoly solution}
\end{center}
\end{table}

\begin{table}
\begin{center}
\begin{tabular}{l|l}
$\int t^k\mathrm{d}\mu$ & $\tfrac{1}{k+1} + (1-\varepsilon)0^k$\\
$\int y^k\mathrm{d}\mu$ & $\tfrac{(1 - \varepsilon)^{k + 1}}{(k + 1)} - \tfrac{(1 - \varepsilon)^{k + 1} - 1}{\varepsilon(k + 1)} $\\
$\int r^k\mathrm{d}\mu$ & $0^k+(1-\varepsilon)$\\
$\int w^k\mathrm{d}\mu$ & $\varepsilon^k+(1-\varepsilon)$\\
\end{tabular}
\caption{Analytic expressions of the moments.}
\label{tab: analytic moments}
\end{center}
\end{table}

\section{\new{Conclusion}}\label{sec:conclu}

\new{
In this paper we have described a unified methodology to cope with limit phenomena typical of optimal control, namely oscillations, concentrations and discontinuities. Our approach relies numerically on the Lasserre hierarchy of semidefinite programming relaxations, which allows for the application of off-the-shelf computer software and hence opens the possibility for engineering applications. The key mathematical tool are anisotropic parametrized measures, an extension of DiPerna-Majda measures, themselves an extension of Young measures, objects familiar to PDE analysts.

From the numerical solution and the nature of the measure computed, we can deduce whether an integrable optimal control function exists or not:
\begin{itemize}
\item if the control measure is concentrated on the graph of a function, then there exists an integrable optimal control law, and hence there is no oscillation of the control;
\item if the time measure is absolutely continuous w.r.t. the Lebesgue measure, then there is no concentration of the control, and hence no discontinuity of the trajectory.
\end{itemize}

Our approach is global, and hence closer in spirit to the Hamilton-Jacobi-Bellmann  approach rather than the Pontryagin Maximum Principle. We are not aware of first-order optimality conditions for optimal control problems at the level of generality considered in our paper. 

Beyond providing a numerical solution to optimal control problems that are potentially troublesome for alternative numerical methods, we believe that our work can pave the way for the application of the Lasserre hierarchy to other problems of calculus of variations and optimal control, especially subject to PDE constraints. Indeed, we would like to emphasize that the theory of anisotropic parametrized measures can also deal with vector-valued multi-dimensional problems.
}


\section*{Acknowledgements}

This work was supported by a CNRS PICS project.
The research of M. Kru\v z\'{\i}k was also supported by the GA\v CR grant 17-04301S.
The research of D. Henrion and T. Weisser was also supported by the ERC advanced grant Taming.

\bibliographystyle{acm}

\end{document}